\documentclass[a4paper,10pt]{article}

\def\eps{\varepsilon}

\def\N{\mathbf N}
\def\R{\mathbf R}
\usepackage[T2A]{fontenc}
\usepackage[cp1251]{inputenc}
\usepackage[english,russian]{babel}
\usepackage[tbtags]{amsmath}
\usepackage{amsfonts,amssymb}

\usepackage{mathrsfs}

\usepackage{graphicx}

\date{}
\begin{document}
\Large

\title{ Тензорно простой спектр унитарных потоков}
\author{Рыжиков В.В.}

\maketitle

\begin{abstract}

В заметке предложены унитарные потоки $T_t$ динамического происхождения такие, что  для всякого  счетного подмножества 
$Q\subset (0,+\infty)$ тензорное произведение  $\bigotimes_{q\in Q} T_q $ имеет однократный спектр.
Типичные потоки, сохраняющие  сигма-конечную меру, обладают этим свойством.

Библиография: 6 названий, УДК: 517.987, \ MSC: Primary 28Y05; Secondary 58F11

Ключевые слова и фразы: \it унитарный поток, слабое замыкание, спектр, тензорное произведение операторов.
\rm

\end{abstract}

%%%%%%%%%%%%%%%%%%%%%%%%%%%%%%%%%%%%%%%%%%%%%%%%%%%%%%%%%%%%%%%%%%%%%%%%%%%%%%%%%%%%%%%

\section{Введение }

В \cite{A}   показано, что для типичного автоморфизма $T$  вероятностного пространства   спектр 
 произведения $T\otimes T^2\otimes T^3\otimes\dots$ простой. 
Этот результат стимулировал поиск унитарного потока с аналогичным (но более   тонким)  спектральным свойством.
В  \cite{LR}  были предъявлены сохраняющие меру  потоки $T_t$  на вероятностном пространстве такие,  что 
для всех $a>1$ произведение  $T_t\otimes T_{at}$ обладало  простым спектром, что обеспечивалось спецификой
 слабого замыкания  потока $T_t$.   Слегка модифицируя эти примеры, мы покажем
существование унитарного потока   $T_t$, для которого спектр  произведений  
$\bigotimes_{q\in Q} T_q $   простой  для всякого   конечного и, следовательно, для всякого счетного множества 
$Q\subset (0,+\infty)$.  Спектр такого потока будем называть тензорно простым.
 
 Потоком $T_t, t\in \textbf{R}$ на пространстве Лебега $(X, \mu)$
 называется непрерывное вложение  $\textbf{R}$ в группу   всех 
автоморфизмов пространства $(X, \mu)$.
Поток $T_t$, сохраняющий меру, индуцирует унитарный поток  на гильбертовом пространстве $L_2(X, \mu)$:
  $T_tf(x)= f(T_tx)$ (оба потока обозначаем одинаково).
Поток $T_t$ имеет простой спектр, если существует вектор $f\in L_2(X, \mu)$ такой, 
что линейное и топологическое замыкание множества  $\{T_tf:t\in \textbf{R}\}$ есть $L_2(X, \mu)$.

\bf Лакунарная жесткость. \rm  Унитарный поток $T_t$  обладает  свойством лакунарной жесткости, если 
  для некоторой последовательности 
$r_j \to \infty$ имеет место сходимость 
$T_{r_j} \to I,$
причем 
$
T_{a_j} \to 0,
$
для  всякой   последовательности 
$\left\{ a_j \right\}$, $a_j \in \left[\eps r_j, (1-\eps)r_j \right]$,
при   $\eps \in (0, \frac 1 2)$.

 В этом случае для произвольного $\alpha>0$ найдется последовательность  $n_j\in\N$ такая,  что 
числа ${\alpha n_j}-r_j$ ограничены и для некоторого $u\in\R$ выполнено 
$$T_{\alpha n_j}\to T_u,$$ 
 $$T_{a_j} \to 0, \ a_j \in \left[\eps \alpha n_j, (1-\eps)\alpha n_j \right].$$  
Такую последовательность  $\alpha n_j$  будем называть  для краткости \it лакунарной. \rm

\bf Специальный слабый предел. \rm Оператор $P$ является специальным 
слабым пределом для потока $T_t$,
если для любого конечного набора положительных  чисел $ \alpha_1,\dots,\alpha_m$ найдется последовательность 
$ n_{j}\in\N$ такая, что  для всех $k=1,2,\dots,m$ имеет место слабая операторная сходимость
$$T_{\alpha_k n_{j}} \to  P.$$
В дальнейшем существенно используется тот факт, что специальные слабые пределы  образуют полугруппу.
Обозначим  $ {\cal P}(T_\beta)= \frac {T_{\beta}+{2} I + T_{-\beta}}{4}$.

\vspace{3mm}
\bf Теорема. \it Пусть  унитарный поток $T_{ t}$ с простым непрерывным спектром  
обладает свойством лакунарной жесткости.
Если все операторы   $ {\cal P}(T_\beta)$, $\beta\in \R$, являются специальными
слабыми пределами для потока $T_{ t}$, то при \\ $0<\alpha_1 < \dots <\alpha_n$  оператор 
 $T_{\alpha_1 }\otimes T_{\alpha_2 }\otimes \dots \otimes T_{\alpha_n }$   
  имеет простой спектр.\rm

 \section{Доказательство теоремы} 
 Пусть ${3\alpha n_j}$ -- лакунарная  последовательность и $T_{3\alpha n_j}\to T_u.$ 
Напомним, что для   последовательности $\{a_j\}$, удовлетворяющей условию
 $a_j \in \left[\eps 3\alpha n_j, (1-\eps)3\alpha n_j \right]$ при  $\eps \in (0, \frac 1 2)$, 
 операторы $T_{a_j}$ слабо сходятся к  $0$.  Следовательно,   $T_{b n_j} \to 0$  при $0<|b|<3\alpha$,
чем мы воспользуемся ниже.

В дальнейшем мы будем обозначать через $W(U)$ алгебру фон Неймана, порожденную оператором $U$. Она является наименьшей 
алгеброй, содержащей $U$, инвариантной относительно эрмитова сопряжения и замкнутой в слабой операторной топологии.
  
Так как  $T_{2\alpha n_j}{\cal P}(T_{\alpha n_j +s})\in W(T_\alpha)$  и для любого $s\in \R$ выполняется
$$ T_{2\alpha n_j}\frac 1 {4}\left(T_{-\alpha n_j-s} +2I+T_{\alpha n_j+s}\right)\to 0+0+ \frac 1 {4}T_{u+s},$$
получаем $\frac 1 {4}T_{t_1}\in W(T_\alpha)$ для всех $t_1\in \R$ ($t_1=u+s$).
Тогда  циклический вектор $f$ для потока $\{T_t\}$ является циклическим вектором для оператора $T_\alpha$.
Значит, для всех $\alpha\neq 0$ спектр оператора $T_\alpha$ простой.

Убедимся в том, что при $0< \alpha_1 <\alpha_2 $   спектр  произведения 
 $T_{\alpha_1 }\otimes T_{\alpha_2 }$ на  $H^{\otimes 2}$  однократен. 
Рассмотрим лакурнарную последовательность
$T_{(\alpha_1+\alpha_2) n_j}\to T_u$, $u>0$, и последовательность операторов
$$Q_j=T_{\alpha_1 n_j}{\cal  P}(T_{\alpha_1 n_j}) \otimes                  
T_{\alpha_2 n_j}{\cal  P}(T_{\alpha_1 n_j}). $$                   
  Так как 
$$T_{(\alpha_1-\alpha_2) n_j}, \ T_{\alpha_1 n_j}, \ T_{2\alpha_1 n_j},\  , 
T_{(\alpha_2-\alpha_1) n_j},  T_{\alpha_2 n_j} \ \to 0,$$
получим $$Q_j\to \frac {  I\otimes T_{u}}{4^2}.$$ 
Учитывая, что   $T_{t}\in W(T_{u })$ для всех $t\in\R$, последовательно устанавливаем принадлежность
 алгебре $W(T_{\alpha_1 }\otimes T_{\alpha_2 })$ операторов
$I\otimes T_{u},\ I\otimes T_{t_2}$,  $T_{\alpha_1}\otimes I,\ T_{t_1}\otimes I,$ 
 $T_{t_1}\otimes T_{t_2}, \ t_1, t_2\in\R.$

Пусть    $f$  -- циклический вектор  потока $T_t$, действующего на пространстве $H$, 
тогда   пространство, инвариантное относительно  оператора $T_{\alpha_1 }\otimes T_{\alpha_2 }$ и содержащее
вектор ${f\otimes f}$,
содержит все векторы вида  $T_{t_1 }f\otimes T_{t_2  }f$, $\ t_1, t_2\in\R$, следовательно,
оно совпадет с $H{\otimes } H$.  Это доказывает, что спектр оператора $T_{\alpha_1 }\otimes T_{\alpha_2 }$ простой.

 Установим  однократность спектра  оператора 
$T_{\alpha_1 }\otimes T_{\alpha_2 }\otimes T_{\alpha_3 }$
%в пространстве $H^{\otimes 3}$ 
при выполнении условия  $0< \alpha_1 <\alpha_2 <\alpha_3$.
Рассмотрим лакурнарную последовательность
$T_{(\alpha_1+\alpha_2+\alpha_3) n_j}\to T_u$, $u>0$ и последовательность операторов
$$Q_j=T_{\alpha_1 n_j}{\cal  P}(T_{\alpha_1 n_j}) {\cal  P}(T_{\alpha_2 n_j})\otimes      
          T_{\alpha_2 n_j}{\cal  P}(T_{\alpha_1 n_j}) {\cal  P}(T_{\alpha_2 n_j})\otimes
T_{\alpha_3 n_j}{\cal  P}(T_{\alpha_1 n_j}) {\cal  P}(T_{\alpha_2 n_j}). $$  
Каждый сомножитель в таком тензорном произведении является суммой девяти  операторов,
из которых    некоторые равны $I$, а остальные имеют вид $T_{\beta n_j}$, $\beta\neq 0$. Все $T_{\beta n_j}$, 
кроме   $T_{(\alpha_1+\alpha_2+\alpha_3) n_j}$, стремятся слабо к 0 в силу лакунарности последовательности 
$(\alpha_1+\alpha_2+\alpha_3) n_j$, так как $|\beta|<\alpha_1+\alpha_2+\alpha_3$.

  Если $ \alpha_1+\alpha_2-\alpha_3\neq 0$, получаем
$$Q_j\to c I\otimes\ I\otimes T_u, \ c>0,$$
а в случае соотношения    $ \alpha_1+\alpha_2=\alpha_3$  имеем   
$$Q_j\to bI\otimes I\otimes I + c I\otimes I\otimes T_u, \ b,c>0.$$ 
Теперь последовательно получаем
 принадлежность алгебре $ W(T_{\alpha_1 }\otimes T_{\alpha_2 }\otimes T_{\alpha_3 })$  операторов
$$I\otimes I\otimes T_{u},\ \
I\otimes I\otimes T_{t_3}, \ \
T_{\alpha_1-\alpha_3}\otimes T_{\alpha_2-\alpha_3}\otimes I, \ \
T_{t_1}\otimes T_{t_2}\otimes I, \
T_{t_1}\otimes T_{t_2}\otimes T_{t_3}$$
для всех $t_1, t_2,t_3\in\R$.
Из сказанного вытекает, что вектор  ${f\otimes f\otimes f}$ является циклическим для оператора 
$T_{\alpha_1 }\otimes T_{\alpha_2 }\otimes T_{\alpha_3 }$, следовательно, спектр этого оператора  простой. 
 
Случай $n>3$ принципиально не отличется от случая $n=3$. 
Пусть 
 для всевозможных  наборов различных положительных $\tilde\alpha_1, \dots,\tilde\alpha_{n-1}$
установлено,  что  для всех $t_1,\dots,t_{n-1}\in\R$ выполнено
$$\bigotimes_{k=1}^{n-1} T_{t_k}   \in 
W\left(T_{\tilde\alpha_1}\otimes \dots \otimes  T_{\tilde\alpha_{n-1}}\right).$$

Пусть $0<\alpha_1 < \dots <\alpha_n$, рассмотрим лакурнарную последовательность
$T_{(\alpha_1+\dots+\alpha_n) n_j}\to T_u$, $u>0$, и определим операторы
$$Q_j=\bigotimes_{k=1}^n \left(T_{\alpha_k n_j}\prod_{m=1}^{n-1}{\cal P}(T_{\alpha_m n_j })\right).$$
Рассуждения, аналогичные предыдущим, приводят к слабой сходимости вида  
$$Q_j\to b_n I\otimes \dots \otimes I + c_n I\otimes \dots \otimes I\otimes T_u, \ b_n\geq 0,c_n>0.$$
Случай $b_n>0$ возникает при наличии соотношений вида 
$$\sum_{i=1}^{n-1}d_i\alpha_i=\alpha_n, \  d_i\in\{-1,0,1\}. \eqno (\ast)$$
Полученные  слабые пределы и индуктивное предположение  обеспечивают  принадлежность  операторов 
$$\left(\bigotimes_{k=1}^{n-1} I\right)\otimes T_{t_n},  \ \  
\left(\bigotimes_{k=1}^{n-1} T_{t_k} \right)  \otimes I, \ \ \bigotimes_{k=1}^{n} T_{t_k}$$
к алгебре $W\left(\bigotimes_{k=1}^n T_{\alpha_k}\right)$ для всех $t_1,\dots,t_n\in\R$.   
Отсюда вытекает однократность  спектра оператора $\bigotimes_{k=1}^n T_{\alpha_k}$. Теорема доказана.

\section{ Заключительные замечания}
\bf Примеры. \rm   Потоки, обладающие  обсуждаемыми  слабыми пределами  предъявлены  в \cite{LR}. Сделаем небольшое уточнение: параметры $r_j$ в цитированной статье удовлетворяли условию 
$r_j = 2^{n_j} - 1 > h_j^3$.  Теперь  наложим  ограничение
 $r_j = 2^{n_j} - 1 > h_j^j$,
что обеспечит для всех $n$ (а не только для $n=2$)  слабые сходимости вида
$$T_{\alpha_1 t_j}\otimes \dots \otimes T_{\alpha_n t_j} \ \to  \  {\cal  P}(T_{\beta})^{\otimes n}.$$
Сохраняя свойства слабого замыкания, указанные примеры  
 модифицируем в потоки с тензорно простым спектром, у которых   фазовое пространство  имеет бесконечную меру. Для этого    все параметры $s_j (i)$,  $i=1,2,\dots,r_j,$, конструкций из \cite{LR} увеличим на $h_j$. 

\vspace{2mm} 
\bf Типичность. \rm Потоки, сохраняющие сигма-конечную меру $\mu$, оснащаются полной метрикой $d$:
$$  d(\{ R_t\}, \{ T_t\})= \max \{\rho (R_s, T_s)\, :\, s\in [0,1]\},$$
где $$\rho (R_s, T_s)=\sum_i \frac { \mu( R_sA_i\,\Delta\, T_sA_i)+\mu( R_{-s}A_i\,\Delta\, T_{-s}A_i)}{2^i}$$
для фиксированного набора множеств  $A_i$, плотного в семействе всех множеств конечной меры и удовлетворяющего условию $\sum_i { \mu( A_i)}2^{-i}<\infty$.
Множество потоков, удовлетворяющих условиям теоремы, содержит  плотное $G_\delta$-множество
 относительно метрики $d$. Несложное доказательство этого факта
использует перенос  рассуждений из  \cite{K} на наш случай. Таким образом,  \it потоки с тензорно простым спектром являются типичными в пространстве потоков, сохраняющих   сигма-конечную меру. \rm

Свойство тензорной простоты спектра сохраняющего меру  действия локально компактной коммутативной группы без кручения  определим так: простой спектр имеет тензорное произведение всякого конечного 
набора операторов, входящих в действие, среди которых нет одинаковых, взаимно обратных и  тождественного. Верно ли, что типичные действия таких групп обладают тензорно простым спектром?

Типичные автоморфизмы вероятностного пространства включаются в многомерные потоки \cite{T}.
Верно ли, что последние  обладают тензорно простым спектром? О разнообразных свойствах типичных
автоморфизмов см. недавнюю работу \cite{GTW}.  

\vspace{2mm} 
\bf Гауссовские и пуассоновские потоки. \rm В связи с  изучением спектральных свойств  гауссовских потоков и пуассоновских надстроек (см.\cite{KSF})   представляет  интерес  следующая задача.  \it Существует ли  поток $T_t$, сохраняющий сигма-конечную меру,  для которого  поток  $\exp(T_t)$ обладает тензорно простым спектром? \rm
Напомним, что 
$$\exp(T_t)= \bigoplus_{n = 0}^{\infty} T_t^{\odot n},$$
где $ T^{\odot 0}$ -- одномерный тождественный оператор,
$ T^{\odot n}$ - симметрическая тензорная  $n$-степень  оператора $T$.
Можно показать, что \it  найдется поток 
$T_t,$ сохраняющий сигма-конечную меру, такой, что
для всякого набора рационально независимых чисел $\alpha_1, \dots,\alpha_n \,\in\R$   
 тензорное  произведение 
$$G=\exp(T_{\alpha_1 })\otimes \dots \otimes \exp(T_{\alpha_n })$$
     имеет простой спектр.\rm

Поясним, как используется рациональная независимость чисел $\alpha_i>0$.  
Пусть $0<\alpha_1 < \dots <\alpha_n$.
Cоотношения вида $(\ast)$ для $\alpha_i$ не возникают, поэтому, как  отмечалось выше, 
слабое замыкание степеней  произведения
$T_{\alpha_1 }\otimes \dots \otimes T_{\alpha_n }$ содержит 
операторы вида $cI\otimes \dots \otimes I\otimes T_{t_n}$. По индукции  показываем, что 
для всевозможных  $t_1,\dots, t_n\in\R$  ненулевые операторы
вида $a T_{t_1 }\otimes  \dots \otimes T_{t_n }$ также лежат в упомянутом 
слабом замыкании. При помощи  таких слабых пределов несложно 
установить однократность  спектра  и 
дизъюнктность всех различных произведений вида
$$T_{\alpha_1 }^{\odot m_1}\otimes \dots \otimes T_{\alpha_n }^{\odot m_n}, \ \ m_i\geq 0. $$
Из этого  непосредственно следует, что спектр произведения $G$ простой. 

Анализ общего случая (без условия рациональной независимости чисел $\alpha_i$) приводит к следующей  задаче: 
пусть $ U, V, W$ -- унитарные операторы с непрерывным спектром такие, что произведение 
$U\otimes V\otimes W$ имеет простой спектр. Какие дополнительные условия
на спектры операторов $ U, V, W$ гарантируют однократность спектра оператора 
$ (U\otimes V)\oplus (V\otimes W) $?


\normalsize
\begin{thebibliography}{99}
    \bibitem{A} O. Ageev, The homogeneous spectrum problem in ergodic theory, Invent. Math., 160:2 (2005), 417--446 

\bibitem{GTW} E.~Glasner, J.-P.~Thouvenot, B.~Weiss. On some generic classes of ergodic measure preserving transformations. Trans. Moscow Math. Soc. 82 (2021), 15--36

\bibitem{K} И.~В.~Климов. Простой спектр тензорных произведений и
типичные свойства сохраняющих меру потоков.
Матем. заметки. 104:6(2018), 942--944; 
Simple Spectrum of Tensor Products and Typical Properties of Measure-Preserving Flows, Math. Notes, 104:6 (2018), 927--929

\bibitem{KSF} И.~П.~Корнфельд, Я. Г. Синай, С. В. Фомин, Эргодическая теория, Наука, М., 1980

\bibitem{LR}М.~С.~Лобанов, В.~В.~Рыжиков.
Специальные слабые пределы и простой спектр тензорных произведений потоков,
Матем. сб. 209:5(2018), 62--73; \ 
Special weak limits and simple spectrum of the tensor products of flows. Sb. Math., 209:5 (2018), 660--671


\bibitem{T} С.~В.~Тихонов. Вложения действий решетки в потоки с многомерным временем, Матем. сб., 197:1 (2006), 97--132; Embedding lattice actions in flows with multidimensional time, Sb. Math., 197:1 (2006), 95--126  

\end{thebibliography}
\end{document}